\documentclass[11pt,reqno]{amsart}
\usepackage{amsmath,amssymb}
 \makeatletter
 \evensidemargin  \oddsidemargin
 \makeatother

\begin{document}
\thispagestyle{empty}
 \setcounter{page}{1}

\begin{center}
{\large\bf  AN OPERATOR PRODUCT INEQUALITIES FOR POLYNOMIALS}

\vskip.20in

{\bf M. Ahmadi Baseri, M. Bidkham and M.Eshaghi} \\[2mm]

\address{Department of Mathematics,
Semnan University,\\ P. O. Box 35195-363, Semnan, Iran}

\email{m\underline ~~ ahmadi\underline ~~  baseri@yahoo.com,
mbidkham@walla.com,madjid.eshaghi@gmail}

\vskip 5mm

\end{center}
\vskip 5mm \noindent {\footnotesize{\bf ABSTRACT.}}
 Let $P(z)$ be a polynomial of degree $n\geq 1$. In this paper we define an operator
$B$, as
 following,
$$B[P(z)]:=\lambda_0 P(z)+\lambda_1 (\frac{nz}{2}) \frac{P'(z)}{1!}+\lambda_2 (\frac{nz}{2})^2 \frac{P''(z)}{2!},$$
where $\lambda_0,\lambda_1$ and $\lambda_2$ are such that all the
zeros of
$$u(z)=\lambda_0 +c(n,1)\lambda_1 z+c(n,2) \lambda_2 z^2$$
lie in half plane
$$|z|\leq |z-\frac{n}{2}|$$
and obtain a new generalization of some well-known results.

 \vskip.10in

  \newtheorem{df}{Definition}[section]
  \newtheorem{rk}[df]{Remark}
   \newtheorem{lem}[df]{Lemma}
   \newtheorem{thm}[df]{Theorem}
   \newtheorem{pro}[df]{Proposition}
   \newtheorem{cor}[df]{Corollary}
   \newtheorem{ex}[df]{Example}

 \setcounter{section}{0}
 \numberwithin{equation}{section}

\vskip .2in

\begin{center}
\section{Introduction and statement of results}
\end{center}
Let $P_n$ be the class of polynomials of degree at most $n$ then,
$$\max_{|z|=1} |P'(z)|\leq n \max_{|z|=1} |P(z)| \eqno
(1.1)$$ and
$$\max_{|z|=R} |{P}(z)|\leq R^n \max_{|z|=1} |P(z)|,\quad R>1. \eqno
(1.2)$$ Inequality (1.1) is an immediate consequence of
Bernstein's theorem on the derivative of a polynomial (see[5]).
Inequality (1.2) is a simple deduction from the maximum modulus
principle (see[11]). If we restrict ourselves to the class of
polynomials having no zeros in $|z|<1$, then the inequality (1.1)
and (1.2) can be respectively replaced by following [9,1]
$$\max_{|z|=1} |P'(z)|\leq \frac{n}{2} \max_{|z|=1} |P(z)| \eqno
(1.3)$$ and
$$\max_{|z|=R} |{P}(z)|\leq \frac{R^n +1}{2} \max_{|z|=1} |P(z)|,\quad R>1. \eqno
(1.4)$$ Recently Aziz and Rather [3] have investigated the
dependence of
$$\displaystyle\left|P(Rz)-\alpha P(z)+\beta\displaystyle\left\{\displaystyle\left(\frac{R+1}{2}\right)^n-|\alpha|
\right\}P(z)\right|on ~\max_{|z|=1}|P(z)|$$ for all real and
complex $\alpha,\beta$ with $| \alpha | \leq1$, $| \beta | \leq 1$
and $R\geq1$, in fact, they proved
\paragraph{\bf Theorem A}If P(z) is a polynomial of degree $n$, then for all real or
complex numbers $\alpha , \beta$ with $| \alpha |\leq 1$,
$|\beta|\leq 1$ and $R\geq 1$,
\begin{align*}
   &\displaystyle\left|P(Rz)
-\alpha P(z)+\beta\displaystyle\left\{\displaystyle\left(\frac{R+1}{2}\right)^n-|\alpha| \right\}P(z)\right| \\
   &\leq \displaystyle\left|R^n -\alpha +\beta
\displaystyle\left\{\displaystyle\left( \frac{R+1}{2}\right)^n
-|\alpha |\right\}\right| |z|^n \max_{|z|=1}|P(z)|~For~|z|\geq
1\hspace{1.5cm} (1.5)
\end{align*}
and
\begin{align*}
&\displaystyle\left|P(R z)-\alpha P(z)+\beta\left\{\left(\frac{R+1}{2}\right)^n-|\alpha| \right\}P(z)\right|\\
&+\displaystyle\left| Q(R z)-\alpha Q(z)+\beta\left\{ \left(\frac{R+1}{2}\right)^n-|\alpha |\right\}Q(z)\right|\\
&\leq\left[\left| R^n -\alpha +\beta \left\{\left(
\frac{R+1}{2}\right)^n -|\alpha |\right\}\right|
|z|^n\right.\\
&+\left.\displaystyle\left|1-\alpha +\beta\left \{\left(
\frac{R+1}{2}\right)^n -|\alpha |\right \} \right|\right]
\max_{|z|=1}|P(z)|\hspace{4.3cm} (1.6)
\end{align*}
 for $|z|\geq1$, where
$Q(z)=z^n \overline{P(\frac{1}{\bar{z}})}$. The results are sharp
and equality in (1.5) and (1.6) holds for $P(z)=\lambda z^n$,
$\lambda \neq 0$. For the class of polynomial having no zeros in
$|z|<1$, we have the following result due to Aziz and Rather which
is a generalization of inequality (1.4).
\paragraph{\bf Theorem B}If $P(z)$ is a polynomial of degree $n$ which does not
vanish in $|z|<1$, then for all real and complex number $\alpha
,\beta$ with $|\alpha|\leq 1$ , $|\beta|\leq 1$ and $R\geq 1$,
\begin{align*}
&\displaystyle\left|P(Rz)-\alpha P(z)+\beta
\left\{\left(\frac{R+1}{2}\right)^n-
|\alpha|\right \}P(z)\right|\\
&\leq \frac{1}{2}\left[\left|R^n -\alpha +\beta \left\{\left(
\frac{R+1}{2}\right)^n
-|\alpha |\right\}\right| |z|^n\right.\\
&\left.+\left|1-\alpha +\beta\left \{\left( \frac{R+1}{2}\right)^n
-|\alpha|\right\}\right|\right] \max_{|z|=1}|P(z)|\hspace{4.5
cm}(1.7)
\end{align*}
for $|z|>1$. Equality in (1.7) occurs for $P(z)=z^n +1$.\\ In this
paper, we consider on operator B which carries $P\in P_n$ in to
$$B[P(z)]:=\lambda_0 P(z)+\lambda_1
(\frac{nz}{2})\frac{P'(z)}{1!}+\lambda_2 (\frac{nz}{2})^2
\frac{P''(z)}{2!},\eqno (1.8)$$ where $\lambda_0$, $\lambda_1$ and
$\lambda_2$ are such that all the zeros of
$$u(z)=\lambda_0 +c (n,1)\lambda_1 z+c (n,2)\lambda_2 z^2 \eqno
(1.9)$$ lie in the half plane
$$|z|\leq |z-\frac{n}{2} | \eqno  (1.10)$$
and prove the following generalization of theorem A and B thus as
well of inequalities (1.1) and (1.2).
\paragraph{\bf Theorem 1.}If $P(z)$ is a polynomial of degree $n$, then for real or complex numbers $\alpha
,\beta$ with $|\alpha|\leq 1$ , $|\beta|\leq 1$ and $R\geq 1$,
\begin{align*}
&\displaystyle\left|B[P(Rz)]-\alpha B[P(z)]+\beta \displaystyle\left\{\displaystyle\left(\frac{R+1}{2}\right)^n-|\alpha|\right \}B[P(z)]\right|\\
&\hspace{0.5cm}\leq\displaystyle\left | (R^n-\alpha
)+\beta\displaystyle\left
\{\displaystyle\left(\frac{R+1}{2}\right)^n -| \alpha
|\right\}\right| |B[z^n]| \max_{|z|=1}|P(z)|~ for ~|z|\geq
1.\hspace{0.2cm} (1.11)
\end{align*}
Equality holds in (1.11) for $P(z)= \lambda z^n, \lambda \neq 0$.
\paragraph{\bf Remark 1.}For $\lambda_0=\lambda_2=0$ in (1.11) and note that in this case all the zeros
of $u(z)$ defined by (1.9) lie in (1.10), we get
\begin{align*}
&\displaystyle\left|R P'(R z)-\alpha P'(z)+\beta
\displaystyle\left\{\displaystyle\left(\frac{R+1}{2}\right)^n-|\alpha| \right\}P'(z)\right|\\
&\hspace{0.6cm}\leq n \displaystyle\left| (R^n-\alpha )+\beta
\displaystyle\left\{\displaystyle\left(\frac{R+1}{2}\right)^n -|
\alpha |\right\}\right| |z|^{n-1} \max_{|z|=1}|P(z)|~ for~ |z|\geq
1.\hspace{0.5mm} (1.12)
\end{align*}
Equality holds in (1.12) for $P(z)= \lambda z^n, \lambda \neq 0$.
If we take $\beta=0$ , $\alpha=1$ and dividing the two sides of
(1.12) by $R-1$ and make $R\rightarrow 1$, we get
$$\max_{|z|=r\geq 1}|zP''(z)+P'(z)| \leq n^2r^{n-1}
\max_{|z|=1}|P(z)|.\eqno (1.13)$$ Equality holds in (1.13) for
$P(z)= \lambda z^n, \lambda \neq 0$. for $\alpha = \beta =0$ and
$R=1$, inequality (1.12) gives
$$|P'(z)| \leq n |z|^{n-1} \max_{|z|=1}|P(z)|.\quad for \quad |z|\geq 1.\eqno (1.14)$$
which in particular gives inequality (1.1).\\ \\
\paragraph{\bf Remark 2.}For $\lambda_1=\lambda_2=0$ in theorem 1 reduces to inequality (1.5).
Next as a application of theorem 1, we prove the following theorem
which is a generalization of a results prove by Rahman [10] , Jain
[6], Aziz and Rather [3].
\paragraph{\bf Theorem 2.}If $P(z)$ is a polynomial of degree $n$, then for real or complex number $\alpha
,\beta$ with $|\alpha|\leq 1$ , $|\beta|\leq 1$ and $R\geq 1$,
\begin{align*}
&\displaystyle\left|B[P(Rz)]-\alpha
B[P(z)]+\beta\displaystyle\left\{
\displaystyle\left(\frac{R+1}{2}\right)^n-|\alpha|\right\}B[P(z)]\right|\\
&+\displaystyle\left|B[Q(Rz)]-\alpha B[Q(z)]+\beta \displaystyle\left\{\displaystyle\left(\frac{R+1}{2}\right)^n-|\alpha|\right\}B[Q(z)]\right|\\
& \leq\displaystyle\left[ \displaystyle\left|R^n-\alpha+\beta
\displaystyle\left\{\displaystyle\left(\frac{R+1}{2}\right)^n -|
\alpha |\right\}\right|
|B[z^n]|\right.\\
&+\left.\displaystyle\left|1-\alpha+\beta
\displaystyle\left\{\displaystyle\left(\frac{R+1}{2}\right)^n -|
\alpha |\displaystyle\right\}\displaystyle\right|
|\lambda_0|\displaystyle \right] \max_{|z|=1}|P(z)|\hspace{3.6cm}
(1.15)
\end{align*}
for $|z|\geq 1$ , where $Q(z)=z^n
\overline{P(\frac{1}{\bar{z}}}).$ If we take
$\lambda_0=\lambda_2=\beta=0$ and $\alpha=1$ in (1.15), we obtain
the following result.
\paragraph{\bf Corollary 1.}If $P(z)$ is a
polynomial of degree $n$, then for all real or complex number
$\alpha$ with $|\alpha|\leq 1$ and $R\geq 1$,
\begin{align*}
|R P'(Rz)&-P'(z)|+|R Q'(Rz)-Q'(z)|\\
&\leq n(R^n-1)|z|^{n-1} \max_{|z|=1}|P(z)|\quad for \quad |z|\geq
1.\hspace{2.9cm}(1.16)
\end{align*}
Equality holds in (1.16) for $P(z)= \lambda z^n, \lambda \neq 0$.
Theorem 2 includes a result due to Rahman [10] as a special case
for $\lambda_1 = \lambda_2 =\alpha =\beta =0$ where as inequality
(1.15) reduces to a result due to Jain [6,Theorem 1] for
$\lambda_1 = \lambda_2 =\alpha =0$. For $\lambda_1 = \lambda_2
=0$, inequality (1.15) reduces to inequality (1.6).\\ Lastly, for
class of polynomial having no zeros in $|z|<1$, we prove the
following generalization of theorem B.
\paragraph{\bf Theorem 3.}If $P(z)$ is a polynomial of degree $n$ which does not
vanish in $|z|<1$, then for all real and complex numbers $\alpha
,\beta$ with $|\alpha| \leq1$, $| \beta| \leq1$  and $R \geq 1$,
then
\begin{align*}
&\displaystyle\left|B[P(Rz)]-\alpha
B[P(z)]+\beta\displaystyle\left\{
\displaystyle\left(\frac{R+1}{2}\right)^n-| \alpha |\right\}B[P(z)]\right|\\
&\leq \frac{1}{2} \displaystyle\left[ \displaystyle\left|
R^n-\alpha+\beta\displaystyle\left\{\displaystyle\left(\frac{R+1}{2}\right)^n
-| \alpha |\right\}\right| |B[z^n]|\right.\\
&+\left. \displaystyle \left|1-\alpha+\beta
\displaystyle\left\{\displaystyle\left(\frac{R+1}{2}\right)^n -
|\alpha |\right\}\right| |\lambda_0| \right]
\max_{|z|=1}|P(z)|\quad for \quad |z|\geq 1.\hspace{1.1cm} (1.17)
\end{align*}
Equality holds in (1.17) for $P(z)=z^n +1$. If we take
$\alpha=\beta=0$ in theorem 3, we get the following result.
\paragraph{\bf Corollary 2.}If $P(z)$ is a polynomial of degree $n$ which does not
vanish in $|z|<1$, then for $R\geq 1$,
$$|B[P(Rz)]|\leq \frac{1}{2} \{R^n |B[z^n]|+|\lambda_0|\} \max_{|z|=1}
|P(z)| \quad for \quad |z|\geq 1.\eqno (1.18)$$ The result is
sharp and equality holds for $P(z)=z^n+1$.  For $R=1$, inequality
(1.18) reduces to a results due to Shah and Liman [12].
\paragraph{\bf Remark 3.}Theorem 3 includes some well-known inequalities as
special case. For example inequality (1.17) reduces to a result
due to Aziz and Rather [4] for $\lambda_1 = \lambda_2 =\beta =0$.
For $\lambda_1 = \lambda_2 =\alpha =0$ inequality (1.17) reduces
to result due to Jain [7] where as for $\lambda_1 = \lambda_2 =0$
inequality (1.18) reduces to
$$|P(Rz)|\leq \frac{1}{2} \{R^n +1 \} \max_{|z|=1}
|P(z)| \quad R \geq 1.$$ If we take $\lambda_0 = \lambda_2 =\alpha
=\beta =0$ , inequality (1.17) reduces to inequality (1.3). \vskip
15mm \vskip 5mm
\begin{center}
\section{ Lemmas}
\end{center}
For the proofs of the theorems, we need the following lemmas.
\paragraph{\bf Lemma 1} If $P(z)$ is a polynomial of degree $n$ having all its
zeros in \\$|z| \leq k (k \leq 1)$, then for every $R>1$,
$$|P(Rz)| \geq\displaystyle\left (\frac{R+k}{1+k}\right)^n |P(z)| \quad for \quad
|z|=1.\eqno (2.1)$$ This lemma was proved by Aziz [2]. If we take
$k=1$ and $P(z) \neq 0$ in $|z| \geq 1$ in the above lemma, then
we can prove the following lemma.
\paragraph{\bf Lemma 2} If $P(z)$ is a polynomial of degree $n$ having all its
zeros in $|z| < 1$, then for every $R > 1$
$$|P(Rz)| > \displaystyle\left(\frac{R+1}{2}\right)^n |P(z)| \quad for \quad
|z|=1.\eqno (2.2)$$
\paragraph{\bf Proof .}Let for $|z_i| >1 (i=1,2,... n), F(z)=z^n \overline{P(\frac{1}{\bar{z}})} =c \prod^{n}_{i=1} (z-z_i)$ for $z=\rho e^{i \theta} (0 \leq \rho < 1)$, we have $| \frac{z-z_i}{e^{i \theta}-z_i}| > \frac{\rho+1}{2} (-\pi \leq \theta \leq \pi)$, therefor
$$\displaystyle\left| \frac{F(\rho z)}{F(z)}\right| > \displaystyle\left(\frac{\rho+1}{2}\right)^n, \quad |z|=1$$
finally for $ \rho=\frac{1}{R}$, we get the desired result.
\paragraph{\bf Lemma 3}If all the zeros of a polynomial $P(z)$ of degree $n$ lie
in a circle $ |z| \leq 1$, then all the zeros of the polynomial
$B[P(z)]$ also lie in the circle $ |z| \leq 1$.
\paragraph{\bf Lemma 4.}If $P(z)$ is a polynomial of degree $n$ such that $P(z)\neq 0$, in
$|z| <1$, then
$$|B[P(z)]| \leq |B [Q(z)]| \quad for \quad
|z| \geq1.\eqno (2.3)$$ where $Q(z)=z^n
\overline{P(\frac{1}{\bar{z}})}.$
\paragraph{\bf Lemma 5}If $P(z)$ is a polynomial of degree $n$, then for $ |z| \geq 1$,
$$|B[P(z)]|+|B[Q(z)]| \leq \{ | B [z^n]|+|\lambda_0| \} \max_{|z|=1} |P(z)|. \eqno  (2.4)$$
The above three lemmas are due to Shah and Liman [12].
\paragraph{\bf Lemma 6}If $P(z)$ is a polynomial of degree $n$ which does not vanish
in $|z|< 1$, then for all real or complex numbers $\alpha,\beta$
with $| \alpha |\leq 1, | \beta |\leq 1$ and $R\geq 1$,
\begin{align*}
&\displaystyle\left|B[P(Rz)]-\alpha B[P(z)]+\beta
\displaystyle\left\{\displaystyle\left(\frac{R+1}{2}\right)^n-|
\alpha |
\right\}B[P(z)]\right|\\
&\leq \displaystyle\left|B[Q(Rz)]-\alpha B [Q(z)]+\beta
\displaystyle\left\{\displaystyle\left(\frac{R+1}{2}\right)^n
-|\alpha |\right\}B[Q(z)]\right| \hspace{2.1cm} (2.5)
\end{align*}
for $|z| \geq 1$ where $Q(z)=z^n \overline{P(\frac{1}{\bar{z}})}.$
\paragraph{\bf Proof.}If $P(z)\neq 0$ in $|z|<1$, then by lemma 4 we have $|B[P(z)]| \leq |B[Q(z)]|$ for $|z| \geq 1$ and hence for $R=1$, we have nothing to prove . For $R>1$, since $|P(z)|=|Q(z)|$ for
$|z|=1$, it follows by Rouche's theorem that for every real or
complex number $\lambda$ with $| \lambda |>1$, the polynomial
$T(z)=P(z)-\lambda Q(z)$ does not vanish in $|z|>1$, with at least
one zero in $|z|<1$. Let $T(z)=(z-re^{i\delta})F(z)$ where $r<1$
and $F(z)$ is a polynomial of degree $n-1$ having no zeros in
$|z|>1$. Applying lemma 1 with $k=1$, for every $R>1$, $0\leq
\theta \leq 2\pi$
\begin{align*}
|T(Re^{i \theta})|  &\geq | Re^{i \theta} -re^{i\delta}|\displaystyle\left (\frac{R+1}{2}\right)^{n-1} |F(e^{i \theta})|\\
&=\displaystyle\left(\frac{R+1}{2}\right)^{n-1}\displaystyle\left
| \frac{Re^{i \theta}-re^{i \delta}}{e^{i \theta}-re^{i\delta
}}\right| |(e^{i \theta}-re^{i \delta} )F(e^{i
\theta})|\\
&\geq\displaystyle\left (\frac{R+1}{2}\right)^{n-1}
\displaystyle\left(\frac{R+r}{1+r}\right) |T(e^{i \theta})|
\end{align*}
or\\
$$\displaystyle\left(\frac{r+1}{R+1}\right) |T(Re^{i \theta})| \geq\displaystyle\left (\frac{R+1}{2}\right)^{n-1} |T(e^{i
\theta})|,  R>1  \quad and \quad 0\leq \theta \leq 2\pi ,\quad
(2.6)$$ since $R>1>r$, hence $T(Re^{i \theta})\neq 0$ and
$(\frac{2}{R+1})>(\frac{r+1}{R+r})$, from inequality (2.6), we
have
$$|T(Rz)| >\displaystyle\left(\frac{R+1}{2}\right)^{n} |T(z)|, \quad |z|=1~,~R>1.\eqno(2.7)$$
Hence for every real or complex number with $|\alpha| \leq 1$, we
have
\begin{align*}
|T(Rz)-\alpha T(z)|&\geq | T(Rz)|- | \alpha | |T(z)|\\
&>\displaystyle\left
\{\displaystyle\left(\frac{R+1}{2}\right)^{n}-| \alpha |\right \}
|T(z)| \quad for ~ |z|=1 ~ and ~ R>1.\hspace{0.2cm} (2.8)
\end{align*}
Since $T(Re^{i \theta})\neq 0$ and $(\frac{R+1}{2})^n >1$, hance
from inequality (2.7), we have
$$|T(Re^{i \theta})| >(|T(e^{i \theta})|, \quad for \quad R>1 ~ and~ \quad 0\leq \theta \leq 2\pi $$
or
$$|T(Rz)| >(|T(z)|, \quad for \quad |z|=1~ and~ \quad R>1.$$
Since all the zeros of $T(Rz)$ lie in $|z|<1$, it follows (by
Rouche's theorem for $| \alpha | \leq 1$) that the polynomial
$T(Rz)-\alpha T(z)$ does not vanish in $|z| \geq1$. Hence from
inequality (2.7)(by Rouche
 theorem for $| \beta | \leq 1$), we have the polynomial
$$S(z)=T(Rz)- \alpha T(z) +\beta \displaystyle\left\{\displaystyle\left(\frac{R+1}{2}\right)^n-| \alpha |\right\}
T(z) $$has all its zeros in $|z|<1$. Therefore, by lemma 3, all
the zeros of $B[S(z)]$ lie in $|z|<1$. Replacing $T(z)$ by $P(z)-
\lambda Q(z)$~and since B is liner, it follows that
\begin{align*}
&\displaystyle\left |B[P(Rz)]-\alpha B[P(z)]+\beta
\displaystyle\left\{\displaystyle\left(\frac{R+1}{2}\right)^n-|\alpha|
 \right\}B[P(z)]\right| \\
& \leq \displaystyle\left|B[Q(Rz)]-\alpha B[Q(z)]+\beta
\displaystyle\left\{\displaystyle\left(\frac{R+1}{2}\right)^n-|\alpha|
\right\}B[Q(z)]\right|\hspace{2.1cm} (2.9)
\end{align*}
for $|z| \geq 1$. If this is not true, then there is a point
$z=z_0$ with $|z_0| \geq 1$, such that
\begin{align*}
&\displaystyle\left|B[P(Rz_0)]-\alpha B[P(z_0)]+\beta
 \displaystyle\left\{\displaystyle\left(\frac{R+1}{2}\right)^n-|\alpha|\right \}B[P(z_0)]\right|\\
& > \displaystyle\left|B[Q(Rz_0)]-\alpha B[Q(z_0)]+\beta
\displaystyle\left\{\displaystyle\left(\frac{R+1}{2}\right)^n-|\alpha|
\right \}B[Q(z_0)]\right|
\end{align*}
since all the zeros of $Q(z)$ lie in~~$|z| \leq 1$, hence (some as
$T(z)$) all the zeros of
$$B[Q(Rz)]-\alpha B[Q(z)]+\beta\displaystyle\left \{\displaystyle\left(\frac{R+1}{2}\right)^n-|\alpha|
\right\}B[Q(z)]\quad lie~ in |z|<1$$ for all real or complex
$\alpha , \beta$ with $|\alpha| \leq 1$, $|\beta| \leq 1$
and~$R>1$. Therefore
$$B[Q(Rz_0)]-\alpha B[Q(z_0)]+\beta \displaystyle\left\{\displaystyle\left(\frac{R+1}{2}\right)^n-|\alpha|
\right\}B[Q(z_0)]\neq 0\quad whit \quad |z_0| \geq 1 ~ we~take$$
$$\lambda= \frac{B[P(Rz_0)]-\alpha B[P(z_0)]+\beta \displaystyle\left\{\displaystyle\left(\frac{R+1}{2}\right)^n-|\alpha|
\right\}B[P(z_0)]}{B[Q(Rz_0)]-\alpha B[Q(z_0)]+\beta
 \displaystyle\left\{ \displaystyle\left(\frac{R+1}{2}\right)^n-|\alpha|\right \}B[Q(z_0)]}$$
so that $|\lambda|>1$ and for this value $\lambda$, $B[S(z_0)]=0$
for $|z_0|\geq 1$, which contradicts the fact that all the zeros
of $B[S(z)]$ lie in $|z|<1$. This proves the desired result.
\vskip 15mm \vskip 5mm
\begin{center}
\section{ Proof Of Theorems}
\end{center}
\paragraph{\bf Proof of Theorem 1.}Let $M=\displaystyle\max_{|z|=1}|P(z)|$, then $|P(z)| \leq M|z^n|$ for
$|z|=1$. If $\lambda$ is any real or complex number with $|\lambda
|>1$, then by Rouche's theorem the polynomial $P_1(z)=P(z)-\lambda
Mz^n$~has all its zeros in $|z|<1$, and by lemma 3, all the zeros
of $B[P_1 (z)]$ lie in $|z|<1$, it follows that
$$| B[P(z)]| \leq M |B[z^n]|\quad for |z|\geq 1$$
If this is not true, then there is a point $z=z_0$ with $|z_0|
\geq 1$, such that
$$| B[P(z)]| > M |B[z^n]|$$
since $B[z^n]\neq 0$ for $|z| \geq 1$, we take
$\lambda=\frac{B[P(z_0)]}{MB[z_0^n]}$ so that $|\lambda |>1$ and
for this value of $\lambda,B[P_1(z_0)]=0$ for $|z_0|\geq 1$, which
contradicts the fact that all the zeros of $B[P_1(z)]$ lie in
$|z|<1$. This completes the proof for $R=1$. Now we assume $R>1$.
On applying lemma 2, we get
$$|P_1(Rz)|> \displaystyle\left (\frac{R+1}{2}\right)^n|P_1(z)|,\quad
|z|=1~and~R>1.\eqno (3.1)$$ Since $P_1(Re^{i \theta})\neq 0$,
$0\leq \theta \leq 2\pi$ and $(\frac{R+1}{2})^n>1$~from above
inequality we have
 $|P_1(Re^{i\theta})|>|P_1(e^{i \theta})|$,~$R>1$.\\Equivalently,
$$|P_1(Rz)|> |P(z)|,\quad for~ |z|=1~and~R>1.$$
For every real and complex number $\alpha$ with $|\alpha| \leq1$
and using inequality (3.1), we have
\begin{align*}
|P_1(Rz)-\alpha P_1(z)|&\geq | P_1(Rz)|-|\alpha| |P_1(z)|\\
&>  \displaystyle\left\{
\displaystyle\left(\frac{R+1}{2}\right)^n-|\alpha|
\right\}|P_1(z)|\quad for \quad |z|=1\quad and~R>1.
\end{align*}
Now using the arguments similar to those used in lemma 6, the
theorem follows.
\paragraph{\bf Proof of Theorem 2.}The result is trivial if
$(R=1)$ (lemma 5), so we suppose that $R>1$. If
$M=\displaystyle\max_{|z|=1}|P(z)|,$ then $|P(z)|\leq M$ for
$|z|=1$. Now for $\lambda$ with $| \lambda|>1$, we have the
polynomial $W(z)=P(z)+\lambda M$ has no zeros in $|z|<1$ and on
applying lemma 6, we get for $|z| \geq1$ and $R>1$,
\begin{align*}
& \displaystyle\left| B[P(Rz)]-\alpha B[P(z)]+\beta
\displaystyle\left\{\displaystyle\left(\frac{R+1}{2}\right)^n-|\alpha|\right\}B[P(z)]\right.\\
&+\left.\lambda\displaystyle\left[1-\alpha+\beta
\displaystyle\left\{\displaystyle\left(\frac{R+1}{2}\right)^n-|\alpha|\right\}\right]\lambda_0
M \right|
\\
& \leq\displaystyle\left| B[Q(Rz)]-\alpha B[Q(z)]+\beta
\displaystyle\left\{\displaystyle\left(\frac{R+1}{2}\right)^n-|\alpha| \right\}B[Q(z)]\right.\\
&\left.+\bar{\lambda}\displaystyle\left[R^n-\alpha+\beta\displaystyle\left\{\displaystyle\left(\frac{R+1}{2}
\right)^n-|\alpha|\displaystyle\right\}\right]B[z^n] M \right|
 \hspace{4.3cm}(3.2)
\end{align*}
where $| \alpha|\leq1$ , $|\beta|\leq1$ and $Q(z)=z^n
\overline{P(\frac{1}{\bar{z}})}$. choosing the argument of
$\lambda $, which is possible by (1.11)~such that
\begin{align*}
&\hspace{-0.3cm} \displaystyle\left|B[Q(Rz)]-\alpha B[Q(z)]
+\beta\displaystyle\left\{\left(\frac{R+1}{2}\right)^n-|\alpha|\right\}B[Q(z)]\right. \\
&+\left.\bar{\lambda}\displaystyle\left[ R^n-\alpha +
\beta\left\{\displaystyle\left(\frac{R+1}{2}\right)^n-|\alpha|\right\}\right] MB[z^n]\right|\\
&=|\lambda|\displaystyle\left|R^n-\alpha+
\beta\displaystyle\left\{\displaystyle\left(\frac{R+1}{2}\right)^n-|\alpha|\right\}\right|M|B[z^n]|\\
&-\displaystyle\left|B[Q(Rz)]-\alpha
B[Q(z)]+\beta\displaystyle\left\{\displaystyle\left(\frac{R+1}{2}\right)^n-|\alpha|
\right\}B[Q(z)] \right|
\end{align*}
we get from (3.2)
\begin{align*}
&\displaystyle\left|B[P(Rz)]-\alpha B[P(z)]+\beta\displaystyle\left \{\displaystyle\left(\frac{R+1}{2}\right)^n-|\alpha| \right\}B[P(z)]\right| \\
&+ |\lambda | |\lambda_0|\displaystyle\left|1-\alpha+\beta
\displaystyle\left\{\displaystyle\left(\frac{R+1}{2}\right)^n-|\alpha|\right \}\right|M\\
&\leq | \lambda| \displaystyle\left|R^n-\alpha+\beta
\displaystyle\left\{\displaystyle\left(\frac{R+1}{2}\right)^n-|\alpha|
\right\}\right|M|B[z^n]|\\
&-\displaystyle\left|B[Q(Rz)]-\alpha B[Q(z)]+\beta
\displaystyle\left\{\displaystyle\left(\frac{R+1}{2}\right)^n-|\alpha|
\right\}B[Q(z)]\right|\hspace{2.3cm} (3.3)
\end{align*}
for $|z|\geq 1$,~$|\alpha |\leq 1$,~$|\beta|\leq 1$,~and $R>1$,
making $|\lambda|\rightarrow 1$ in (3.3 ), we get (1.15). This
completes the proof of theorem 2.
\paragraph{\bf Proof of Theorem 3.}By hypothesis $P(z)$ does not
vanish in $ |z|<1$, therefore by lemma 6 we have
\begin{align*}
&\displaystyle\left|B[P(Rz)]-\alpha B[P(z)]+\beta \displaystyle\left\{\displaystyle\left(\frac{R+1}{2}\right)^n-|\alpha|\right \}B[P(z)]\right| \\
&\leq\displaystyle\left |B[Q(Rz)]-\alpha B[Q(z)]+\beta
\displaystyle\left\{\displaystyle\left(\frac{R+1}{2}\right)^n-|\alpha|
\right\}B[Q(z)]\right|
\end{align*}
or
\begin{align*}
&2\displaystyle\left|B[P(Rz)]-\alpha B[P(z)]+\beta
\displaystyle\left\{\displaystyle\left(\frac{R+1}{2}\right)^n-|\alpha|
\right\}B[P(z)]\right|\\
&\leq\displaystyle\left |B[P(Rz)]-\alpha B[P(z)]+\beta
\displaystyle\left\{\displaystyle\left(\frac{R+1}{2}\right)^n-|\alpha|
\right\}B[P(z)]\right|\\
&+\displaystyle\left|B[Q(Rz)]-\alpha B[Q(z)]+\beta
\displaystyle\left\{\displaystyle\left(\frac{R+1}{2}\right)^n-|\alpha|
\right\}B[Q(z)]\right|.\hspace{2.3cm}      (3.4 )
\end{align*}
On combining (3.4 ) with (1.15), the proof of theorem 3 is
complete. 
\end{document}